\documentclass[12pt]{amsart}
\usepackage{amssymb,amscd,epsfig}

\newtheorem*{theorem*}{Theorem}
\newtheorem{theorem}{Theorem}
\newtheorem{lemma}[theorem]{Lemma}
\newtheorem{prop}[theorem]{Proposition}
\newtheorem*{conjecture*}{Conjecture}
\newtheorem*{corollary*}{Corollary}
\numberwithin{equation}{section}
\numberwithin{theorem}{section}

\newcommand{\R}{{\mathbb R}}
\newcommand{\C}{{\mathbb C}}
\newcommand{\Z}{{\mathbb Z}}
\newcommand{\Q}{{\mathbb Q}}
\renewcommand{\o}{\operatorname}
\newcommand{\E}{{\mathcal E}}

\newcommand{\sut}{\mathsf{SU}(2)}

\newcommand{\slt}{\o{SL}(2,\C)}

\newcommand{\Uu}{\mathsf{U}}
\newcommand{\Hom}{\mathsf{Hom}}
\newcommand{\Hfng}{\Hom(F_n,G)}
\newcommand{\Hfngg}{\Hom(F_n,G)/G}

\newcommand{\tr}{\mathsf{tr}}
\newcommand{\Aut}{\mathsf{Aut}}
\newcommand{\Out}{\mathsf{Out}}
\newcommand{\Inn}{\mathsf{Inn}}
\newcommand{\Ou}{\Out(F_n)}
\newcommand{\Aufn}{\Aut(F_n)}

\newcommand{\A}{\mathcal{A}}

\begin{document}
\title[Outer automorphisms of free groups]
{An ergodic action of the outer automorphism group of a free group}
\author{William M. Goldman}
\address{
Department of Mathematics,
University of Maryland, College Park, MD  20742} 
\email{wmg@@math.umd.edu}
\date{\today} 
\thanks{The author gratefully acknowledges support from
National Science Foundation grants DMS-0405605 
and DMS-0103889 }

\begin{abstract}
For $n>2$, the action of the outer automorphism group of the rank $n$
free group $F_n$ on  
$\Hom(F_n,\sut)/\sut$ is ergodic with respect to the Lebesgue measure class.
\end{abstract}
\maketitle
\section*{Introduction}
Let $F_n$ be a free group of rank $n>1$ and let $G$ 
be a compact Lie group.
Then $\Hfng$ admits a natural volume form which is invariant
under $\Aufn$. This volume form descends to a finite measure
on the {\em character variety\/} $\Hfngg$ which is invariant
under $\Ou$. The purpose of this note is to prove:

\begin{theorem*}
Suppose that $G$ is a connected group locally isomorphic to 
a product of copies of $\sut$ and $\Uu(1)$.
If $n>2$, then the $\Ou$-action on $\Hfngg$ is ergodic.
\end{theorem*}

We conjecture that $\Ou$ is ergodic on each connected component of $\Hfngg$ 
for every compact Lie group $G$ and $n>2$.

When $G=\Uu(1)$, then this action is just the action of $\mathsf{GL}(n,\Z)$
on the $n$-torus $\R^n/\Z^n$, which is well known to be ergodic. In fact,
certain cyclic subgroups of $\mathsf{GL}(n,\Z)$ act ergodicly.

The proof relies heavily on \cite{Erg}, both in its outline
and a key result.
When $n=2$, the action is {\em not\/} ergodic, since it
preserves the function
\begin{align*}
\Hfngg &\xrightarrow{\kappa} [-2,2] \\
[\rho] &\longmapsto \tr([X_1,X_2])
\end{align*}
where $X_1,X_2$ are a pair of free generators for $F_2$.
However, for each $-2\le t\le 2$, the action is ergodic on 
$\kappa^{-1}(t)$.

When $\pi$ is the fundamental group of a closed surface, then 
Pickrell and Xia~\cite{PickrellXia} have proved $\Out(\pi)$
is ergodic on $\Hom(\pi,G)/G$ for {\em any\/} compact Lie group $G$.

As in \cite{Erg}, the methods here apply when $G$ is any Lie group
having simple factors $\Uu(1)$ and $\sut$. In particular, since this
class of groups is closed under the operation of taking direct products,
the action of $\Ou$ is ergodic on 
\begin{equation*}
\Hom(F_n,G\times G)/(G\times G) \longleftrightarrow
\Hfngg \times \Hfngg.
\end{equation*}
As in \cite{Erg}, the action of $\Ou$ on $\Hfngg$ is {\em weak-mixing,\/}
that is:
\begin{corollary*}
The only invariant finite-dimensional subrepresentation of the induced
unitary representation of $\Ou$ on $L^2(\Hfngg)$ consists of constants.
\end{corollary*}


I would like to thank David Fisher for pointing out an error in the
original proof of Lemma 3.1 and for many helpful suggestions.

\section{Ergodic theory of the $\sut$-character variety}

Let $\{X_1,\dots,X_n\}$ be a set of free generators for $F_n$ and let
\begin{equation*}
X_0 = X_n^{-1}\dots X_1^{-1}. 
\end{equation*}
 Then $F_n$ is the fundamental group of
an $n+1$-holed sphere $S_{n+1}$, 
where the $X_0, X_1, \dots, X_n$ 
correspond to components of $\partial S_{n+1}$. 
The {\em mapping class group\/} $\Gamma_{n+1}$ of $S_{n+1}$
embeds in $\Ou$ as the subgroup preserving the conjugacy classes
of the cyclic subgroups $\langle X_i\rangle$ for $i=0,\dots, n$.

The proof proceeds as follows. Let 
\begin{equation}\label{eq:function1}
\Hfngg \xrightarrow{f} \R 
\end{equation}
be an $\Ou$-invariant measurable function. We show that $f$
is constant almost everywhere.

The main result of \cite{Erg} applied to the surface $S_{n+1}$ 
gives the following:
\begin{prop}\label{prop:erg}
The mapping
\begin{align*}
\Hfngg &\xrightarrow{t_\partial} [-2,2]^{n+1} \\
[\rho] &\longmapsto 
\bmatrix 
\tr(\rho(X_0)) \\ \vdots \\\tr(\rho(X_{n}))
\endbmatrix
\end{align*}
is an {\em ergodic decomposition\/} for the action of $\Gamma_{n+1}$. 
That is, for every $\Gamma_{n+1}$-invariant measurable function
\begin{equation*}
\Hfngg \xrightarrow{h} \R 
\end{equation*}
there exists a measurable function 
\begin{equation*}
[-2,2]^{n+1}\xrightarrow{H}\R 
\end{equation*}
such that $h = H\circ t_\partial$ almost everywhere.
\end{prop} 

Using the embedding of the  mapping class group
\begin{equation*}
\Gamma_{n+1} \hookrightarrow \Ou
\end{equation*}
as above, the $\Ou$-invariant function $f$ is $\Gamma_{n+1}$-invariant,
and hence factors through $t_\partial$. 

By Proposition~\ref{prop:erg} there exists a function
\begin{equation}\label{eq:function}
[-2,2]^{n+1}\xrightarrow{F}\R.
\end{equation}
such that $f = F \circ t_\partial$, where $f$ is the function
discussed in \eqref{eq:function1}.

%

\section{The case of rank $n=3$}
First consider the case $n=3$. Following the notation of
\cite{BG,Erg}, denote the generators by
\begin{equation*}
A = X_1,\quad B = X_2,\quad C = X_3,\quad D = X_0
\end{equation*}
so that $A,B,C,D$  are subject to the relation 
\begin{equation}\label{eq:relation4}
A B C D = 1. 
\end{equation}
A representation $\rho$ is determined by its values
\begin{equation*}
\big(\rho(A), \rho(B), \rho(C) \big) \in G^3
\end{equation*}
on the generators $A, B, C$ and 
\begin{equation*}
\rho(D)= \rho(C)^{-1} \rho(B)^{-1} \rho(A)^{-1}. 
\end{equation*}

\subsection{Trace coordinates}
The equivalence class $[\rho]$ is determined by the seven functions
\begin{align*}
a & = \tr(\rho(A)) \\ 
b & = \tr(\rho(B)) \\ 
c & = \tr(\rho(C)) \\ 
d & = \tr(\rho(D)) = \tr(\rho(D^{-1}))= \tr(\rho(ABC)) \\ 
x & = \tr(\rho(AB)) \\ 
y & = \tr(\rho(BC)) \\ 
z & = \tr(\rho(CA)) 
\end{align*}
subject to the polynomial relation
\begin{align}\label{eq:fourholes}
x^2 + & y^2 + z^2   + xyz  = \\
& \qquad (ab+cd)x + (ad+bc)y + (ac+bd)z  \notag \\
& \qquad \qquad \qquad  +  (4 - a^2 - b^2 - c^2 - d^2 - abcd). \notag
\end{align}
In other words, the $\slt$-character variety of $F_3$ is the hypersurface
in $\C^7$ defined by \eqref{eq:fourholes}.

When $a,b,c,d\in\R$ the topology of the set of $\R$-points is analyzed
in \cite{BG}. In particular the $\sut$-character variety is the union
over the set $V$ of all 
\begin{equation*}
(a,b,c,d)\in [-2,2]^4 
\end{equation*}
satisfying
\begin{align*}
0 \; \ge \; \Delta(a,b,c,d) &=  \left(
2(a^2 + b^2 + c^2 + d^2) - abcd - 16
\right)^2 \\ & \qquad  - (4-a^2) (4-b^2) (4-c^2) (4-d^2) 
\end{align*}
of compact components of the cubic
surface in $\R^3$ satisfying \eqref{eq:fourholes}.
Here is an alternate description with which it is easier to work.

\subsection{Rank two free groups}
\begin{figure}[ht]
\centerline{\epsfig{figure=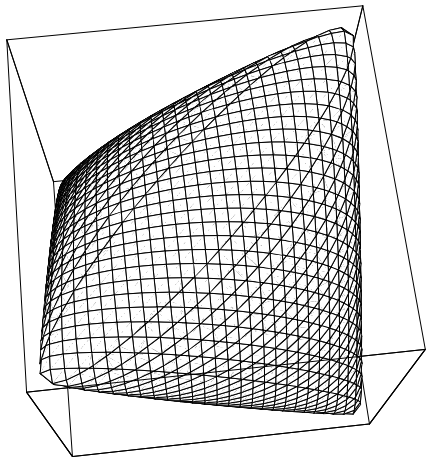,height=6cm}}
\caption{The $\sut$-character variety of a rank two free group is the 
region inside this surface. The surface, a rounded tetrahedron, 
consists of characters of abelian representations. It is a quotient
of the 2-torus by an involution, where the fixed points of the involution
correspond to the vertices of the tetrahedron.}
\label{fig:charv}
\end{figure}

The $\sut$-character variety of $F_2$ is the subset $V_3\subset\R^3$ 
defined by traces
\begin{equation*}
(x_1,x_2,x_3)  \in [-2,2]^3
\end{equation*}
satisfying the inequality
\begin{equation}
x_1^2 + x_2^2 + x_3^2 + x_1x_2x_3  \le 4
\end{equation}
which is depicted in Figure~\ref{fig:charv}.

A quadruple $(a,b,c,d)\in [-2,2]^4$ is the image of an $\sut$-character
if and only if there exists $y\in\R$ such that both triples
$(a,d,y)$ and $(b,c,y)$ lie in $V_3$.

\begin{figure}
\centerline{\epsfig{figure=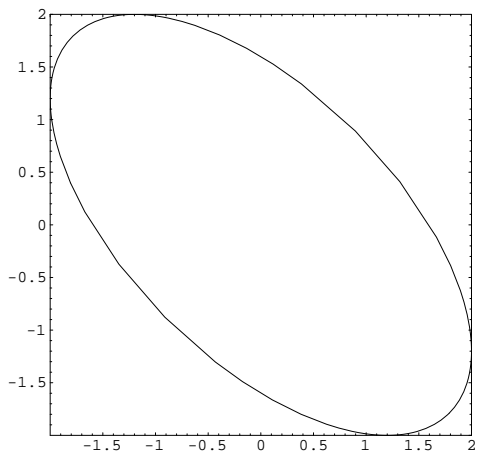,height=6cm}}
\caption{For $y_0 = -1.2$, the set of possible 
$(b,c)\in[-2,2]^2$ for which $(b,c,y_0)$ is the
character of an $\sut$-representation is the interior 
of the ellipse $\E_{1.2}$ inscribed in $b[-2,2]^2$ 
at the four points $(\pm 2, \mp 1.2)$ and $(\mp 1.2,\pm 2)$.}
\label{fig:tellipse}
\end{figure}

\begin{figure}
\centerline{\epsfig{figure=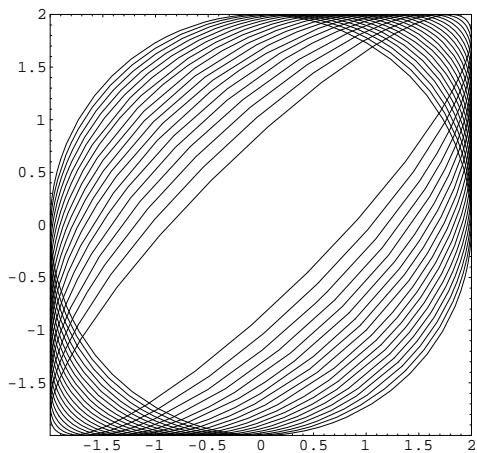,height=6cm}}
\caption{
Here are the ellipses drawn for the interval
$Y = [0,1.8]$. The union of the closed elliptical
regions is the $(b,c)$-projection of the $V_4(a,d)$
where $ad = y_ + y_+ = 1.8$ and $a^2 + d^2=$ (since $y_-=0$),
for example $a=1.69, d = 1.062$.
}
\label{fig:someellipses}
\end{figure}

Using \eqref{eq:fourholes}, we determine the condition that 
$(a,d,y)\in V_3$. Apply \eqref{eq:fourholes}
to $x_1=y,x_2=a,x_3=d$ to see 
that $(a,d,y)\in V_3$ if and only if
$y$ lies in the interval 
\begin{equation*}
Y(a,d) := [y_-(a,d),y_+(a,d)] 
\end{equation*}
with endpoints
\begin{equation*}
y_{\pm}(a,d) := \frac{a d \pm \sqrt{(4-a^2)(4-d^2)} }2 .
\end{equation*}
For any $(a,d)\in[-2,2]^2$, the interval $Y(a,d)$ is nonempty,
so that the restriction of the projection 
\begin{align*}
V_3 &\xrightarrow{\Pi_{a,d}} [-2,2]^2  \\
\bmatrix a \\ b \\ c \\ d \\ x \\ y \\ z \endbmatrix
& \longmapsto \bmatrix a \\ d \endbmatrix
\end{align*}
is onto. Furthermore the fiber $V_3(a,d)$ of the surjection
\begin{align*}
t_\partial(V_3)  & \stackrel{\Pi_{a,d}}\twoheadrightarrow [-2,2]^2 \\
\bmatrix a \\ b \\ c \\ d \endbmatrix & \longmapsto
\bmatrix a \\ d \endbmatrix
\end{align*}
consists of all $(b,c)\in[-2,2]^2$ such that
\begin{equation*}
Y(b,c) \cap Y(a,d) \neq \emptyset.
\end{equation*}
If $-2 < y <2$, then the set of $(b,c)$ such that
$y \in Y(b,c) $ is the closed elliptical region $\bar{\E}_y$
inscribed in the square $[-2,2]^2$ at the four points
\begin{equation*}
(2,y), (y,2), (-2,-y),  (-y,2), 
\end{equation*}
depicted in
Figure~\ref{fig:tellipse}.
If $y=\pm 2$, then the set of $(b,c)$ such that 
$y \in Y(b,c)$ is the line segment 
\begin{equation*}
\big\{ (b,\mp b) \mid -2\le b \le 2 \big\}. 
\end{equation*}
For fixed $a,d\in [-2,2]$, the fiber $\Pi^{-1}\big(V_3(a,d)\big)$ equals
\begin{equation*}
\bigcup_{y\in Y(a,d)} \big\{(b,c) \mid y\in Y(b,c)\big\},
\end{equation*}
depicted in Figure~\ref{fig:someellipses}.


\begin{figure}
\centerline{\epsfig{figure=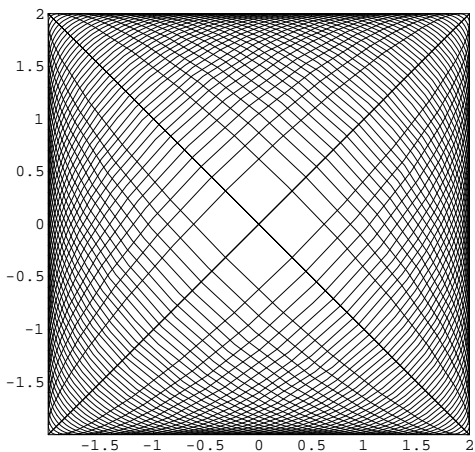,height=6cm}}
\caption{The $\sut$-character variety has three foliations
by ellipses, corresponding to the three coordinate planes.
This figure depicts projections of  the leaves of one of the 
foliations into a coordinate plane. Leaves of the other two
foliations project to to horizontal and vertical line segments,
respectively.
}
\label{fig:charv0}
\end{figure}

\clearpage

\subsection{A non-geometric automorphism}
The automorphism $\alpha\in\Aut(F_3)$ defined by:

\begin{align*}
A & \stackrel{\alpha}\longmapsto A \\ 
B & \stackrel{\alpha}\longmapsto BA^{-1} \\ 
C & \stackrel{\alpha}\longmapsto AC \\ 
D & \stackrel{\alpha}\longmapsto D \\ 
\end{align*}
induces the following automorphism of the character variety:

\begin{equation*}
 \bmatrix a \\ b \\ c \\ d \\ x \\ y \\ z \endbmatrix 
\stackrel{\alpha^*}\longmapsto 
\bmatrix a \\ x \\a c - z \\ d \\ a x - b \\ y \\ c \endbmatrix 
\end{equation*}

Since $\alpha^*$  preserves the coordinates $a,d,y$, this automorphism
restricts to a diffeomorphism on the level sets
$a=a_0, d=d_0$ and $y=y_0$
and the restriction is {\em linear:\/}
\begin{equation*}
\bmatrix x \\ b \\ c \\ z \endbmatrix \stackrel{\alpha^*}\longmapsto
\bmatrix a_0 & -1 & 0 & 0 \\
1 & 0 & 0 & 0 \\
0 & 0 & a_0 & -1  \\
0 & 0 & 1 & 0 \endbmatrix
\bmatrix x \\ b \\ c \\ z \endbmatrix 
\end{equation*}
The following elementary fact (whose proof is omitted) is useful:
\begin{lemma}
Let $-2<a<2$. The linear transformation
\begin{align*}
\R^2 &\xrightarrow{L_a}  \R^2 \\
\bmatrix x \\ y \endbmatrix &\longmapsto
\bmatrix a x - y\\ x \endbmatrix 
\end{align*}
preserves the positive definite quadratic form
\begin{equation*}
\bmatrix x \\ y \endbmatrix \stackrel{Q}\longmapsto  x^2 - a x y + y^2 
\end{equation*}
and lies in the linear flow generated by the vector field
\begin{equation*}
\Upsilon := \bigg(\frac{a}2 x - y \bigg)\partial_x + 
\bigg(x - \frac{a}2 y \bigg)\partial_y.  
\end{equation*}
The trajectories of $\Upsilon$ are the level sets of $Q$, which are ellipses.
$L_a$ is linearly conjugate to a rotation by angle 
\begin{equation}\label{eq:theta}
\theta = 2 \cos^{-1} (a/2). 
\end{equation}
If $\theta\notin \pi\Q$, then $L_a$  has infinite order and 
is ergodic on each $Q$-level set.
\end{lemma}
The $4$-dimensional affine subspaces of $\R^7$ corresponding to levels of
$(a_0,d_0,y_0)$ split as products 
of two affine $2$-planes (corresponding to levels of $(x,b)$ and $(c,z)$
respectively). Evidently the linear map $\alpha^*$ on these
$4$-planes splits as a direct sum of two copies of the linear map
$L_{a_0}$. It preserves the trajectories of the linear vector field
\begin{align*}
\A & = \bigg(\frac{a_0}2 x - b\bigg) \partial_x + 
\bigg(x - \frac{a_0}2 b\bigg) \partial_b \\
& \qquad + \bigg( \frac{a_0}2 c - z\bigg) \partial_c + 
\bigg( c - \frac{a_0}2 z\bigg) \partial_z.
\end{align*}
The zeroes of this vector field consist of the origin 
\begin{equation*}
\bmatrix x \\  b \\ c \\ z \endbmatrix = 
\bmatrix 0 \\  0 \\ 0 \\ 0 \endbmatrix
\end{equation*}
and, when $a_0=\pm 2$, the squares defined by
\begin{equation*}
x = \pm b, \qquad  z = \pm c.
\end{equation*}
All other trajectories are ellipses. 
The transformation $\alpha^*$ acts by rotation along these ellipses through
angle  $\theta$ given by \eqref{eq:theta}.
When $\theta/\pi$ is irrational,
this action is ergodic. Thus, for almost every $a_0\in[-2,2]$, 
the restriction of $\alpha^*$ to
these ellipses is ergodic.  On a set of full measure, the function $f$ is 
constant along the projections of trajectories of $\A$.

Fix $(a_0,d_0)\in (-2,2)^2$ and consider the equivalence
relation $\sim$ on $V_3(a_0,d_0)$ generated by the projections of trajectories
of $\A$. 

\begin{lemma}
For $a_0,d_0\neq \pm 2$, all points in $V_3(a_0,d_0)$ are $\sim$-equivalent.
\end{lemma}
\begin{proof}
Since $V_3(a_0,d_0)$ is connected, it suffices to prove 
that each equivalence class is open.
To this end, 
suppose that $(b_0,c_0)\in V_3(a_0,d_0)$; we show that
every $(b,c)$ sufficiently close to $(b_0,c_0)$ is equivalent to $(b_0,c_0)$.

The image of the tangent vector $\A$ at $(x,b,c,z)$ 
under the differential of the coordinate projection
\begin{equation*}
\bmatrix x \\ b \\ c \\ z \endbmatrix
\stackrel{\Pi_{b,c}}\longmapsto \bmatrix b \\ c \endbmatrix
\end{equation*}
is
\begin{equation*} 
\big(\Pi_{b,c}\big)_*\A = 
\bigg(x - \frac{a_0}2 b\bigg) \partial_b
+ \bigg( \frac{a_0}2 c - z\bigg) \partial_c  
\end{equation*}
The fiber $\Pi_{b,c}^{-1}(b_0,c_0)$ is an interval.
For some (and hence almost every) $(x,z)\in\Pi_{b,c}^{-1}(b_0,c_0)$,
the vector 
\begin{equation*}
\big(\Pi_{b,c}\big)_*\A \; (b_0,x,z,c_0) 
\end{equation*}
is nonzero. Choose such an $(x_0,z_0)$.
For any open neighborhood $U$ of $(x_0,z_0)$, 
the values $\big(\Pi_{b,c}\big)_*\A \; (b_0,x,z,c_0)$
for $(x,z)\in U$, span $\R^2$.

Let $\Phi_t$ denote the flow generated by $\A$ and choose an open
neighborhood $U_0$ of $(b_0,x_0,z_0,c_0)$ in $V_3(a_0,d_0)$.
The differential of the  map
\begin{align*}
\R \times \big(\{b_0\}\times U_0  \times \{c_0\}\big) & \longrightarrow \R^2 \\
(t, b_0,x,z,c_0) &\longmapsto \Pi_{b,c}\big(\Phi_t( b_0,x,z,c_0)\big)
\end{align*}
at $(0,b_0,x_0,z_0,c_0)$ is onto. The inverse function theorem guarantees
an open neighborhood of $(0,b_0,x_0,z_0,c_0)$ mapping onto an open neighborhood
of $(b_0,c_0)$, as desired.
\end{proof}
\noindent
Thus, for almost every $(a_0,d_0)\in[-2,2]^2$,
the function $F$ of \eqref{eq:function} is constant along the level
surfaces $\Pi_{a,d}^{-1}(a_0,d_0)$, and hence factors through the 
projection $\Pi_{a,d}$:
\begin{equation*}
F(a,b,c,d) = F(a,d). 
\end{equation*}
Applying the same argument to the automorphism
\begin{align*}
A & \stackrel{\gamma}\longmapsto CA \\ 
B & \stackrel{\gamma}\longmapsto B \\ 
C & \stackrel{\gamma}\longmapsto C \\ 
D & \stackrel{\gamma}\longmapsto DC^{-1} \\ 
\end{align*}
implies that $F$ factors through the projection $\Pi_{(b,c)}$
and $F$ is almost everywhere constant.

Thus the function $f$, which is invariant under the mapping class
group $\Gamma_4$ and the automorphisms $\alpha_*,\gamma_*$,  must be constant
almost everywhere. This completes the proof of the Theorem when $n=3$.

\section{General rank}

The case of rank $n>3$ follows easily from the case $n=3$. 
The following elementary lemma is useful:
\begin{lemma}\label{lem:ergerga}
Let $G$ be a compact Lie group. Then
$\Ou$ is ergodic on $\Hfngg$ if and only if
$\Aufn$ is ergodic on $\Hfng$.
More generally, let 
\begin{equation*}
\Hfng\xrightarrow{\Pi}\Hfngg
\end{equation*}
denote the quotient map and 
suppose $\Gamma\subset\Aut(F_n)$ is a subgroup 
containing $\Inn(F_n)$.
A measurable $\Gamma$-invariant mapping.
\begin{equation*}
\Hfngg\xrightarrow{f} W 
\end{equation*}
is an ergodic decomposition for $\Gamma$ if and only if
\begin{equation*}
\Hfng\xrightarrow{f\circ\Pi}W
\end{equation*}
is an ergodic decomposition for $\Gamma$.
\end{lemma} 
\begin{proof}
For almost every $\rho\in\Hfng$,
the image $\rho(F_n)$ is dense in $G$. 
Thus $\Inn(F_n)$
is ergodic on each $G$-orbit in $\Hfng$, and the quotient map
$\Pi$ is an ergodic decomposition for the $\Inn(F_n)$-action.
Since $\Inn(F_n)\subset\Gamma$, 
every $\Gamma$-invariant function 
\begin{equation*}
\Hfng \xrightarrow{h} W
\end{equation*}
(where $W$ is a standard measure space)
factors through $\Pi$.
Thus composition with $\Pi$ induces an isomorphism
between the set of equivalence classes of $\Gamma$-invariant measurable maps 
on $\Hfng$ and on $\Hfngg$. In particular the ergodic decompositons
for the respective $\Gamma$-actions are $\Pi$-related, as desired.
\end{proof}
In the following lemma, let $\iota_j$ (for $j=1,\dots,n$) denote the
monomorphism defined by 
\begin{align*}
F_{n-1} &\stackrel{\iota_j}\hookrightarrow F_n \\
X_i &\longmapsto \begin{cases}
X_i \quad\text{~if~} i < j \\ 
X_{i+1} \text{~if~} i \ge j \end{cases}
\end{align*}
for $i=1,\dots,n-1$.
Denote by 
\begin{equation*}
\Aut(F_{n-1)} \stackrel{I_j}\hookrightarrow \Aut(F_n) \\ 
\end{equation*}
the homomorphism defined by 
\begin{equation*}
X_i \stackrel{I_j(\phi)}\longmapsto \begin{cases}
X_j \quad\text{~\qquad\qquad\quad if~} i = j \\ 
\iota_j(\phi(\iota_j^{-1}(X_i))) \text{~\quad if~} i \neq j \end{cases}
\end{equation*}
for $\phi\in\Aut(F_n)$. Clearly $\iota_j$ is equivariant with respect
to $I_j$.


\begin{lemma}\label{lem:erginduct}
Let $n>3$. Suppose that $\Aut(F_{n-1})$
is ergodic on $\Hom(F_{n-1},G)$. 
Let $1\le j\le n$. 
Then 
\begin{align*}
\Hfng & \xrightarrow{t_j} [-2,2]  \\
\rho &\longmapsto \tr\big(\rho(X_j)\big)
\end{align*}
is an ergodic decomposition for the action of 
$I_j(\Aut(F_{n-1}))\times \Inn(G)$.
\end{lemma}
\begin{proof}
Clearly $t_j$ is $I_j\big(\Aut(F_{n-1}\big)$-invariant. 
Let $-2<\tau<2$. The fiber $t_j^{-1}(\tau)$ 
identifies with the set of equivalence classes of
\begin{equation*}
\big((\rho(X_1),\dots,\rho(X_{n-1}),\rho(X_n)\big)  
\end{equation*}
where $\tr\big(\rho(X_j)\big)=\tau$,
that is, the set of $\Inn(G)$-orbits of
\begin{equation*}
\Hom(F_{n-1},G) \times \tr^{-1}(\tau).
\end{equation*}
By Lemma~\ref{lem:ergerga} and the ergodicity
hypothesis, $\Aut(F_{n-1})$ is ergodic on
$\Hom(F_{n-1},G)$. Transitivity of the $\Inn(G)$-action on 
$\tr^{-1}(\tau)$ implies ergodicity of 
$\Aut(F_{n-1})\times\Inn(G)$ acting on 
the subset of $\Hfng$ 
corresponding to $t_j^{-1}(\tau)$, which 
is equivalent to ergodicity of 
$\Out(F_{n-1})$ on $t_j{-1}(\tau)$.
\end{proof}
\noindent
Suppose inductively that $\Aut(F_{n-1})$ is ergodic on
$\Hom(F_{n-1},G)$. By Lemma~\ref{lem:erginduct}, 
any $I_j(\Aut(F_{n-1}))$-invariant measurable function $f$ factors through 
$t_j$ for each $j=1,\dots,n$. We need only consider $j=1$ and $j=n$.
Since 
\begin{equation*}
\Hfngg \xrightarrow{(t_1,t_n)} [-2,2]^2
\end{equation*}
is onto, any two points $\rho,\rho'$ in a full measure subset of
$\Hfng$ can be joined by a third $\rho''$ such that
\begin{align*}
t_1(\rho'') & = t_1(\rho), \\   
t_n(\rho'') & = t_n(\rho').    
\end{align*}
Let $f$ be an $\Aut(F_n)$-invariant measurable function on $\Hfngg$.
Factorization of $f$ through $t_1$ implies that 
\begin{equation*}
f(\rho) = f(\rho'')  
\end{equation*}
and factorization of $f$ through $t_n$ implies that 
\begin{equation*}
f(\rho'') = f(\rho').  
\end{equation*}
Thus $f(\rho)=f(\rho'')$, and $f$ is constant almost everywhere, 
as desired. \qed

\section{Other Lie groups}

The extension of the proof to products of $\sut$ and $\Uu(1)$ proceeds
exactly as in \cite{Erg}. Away from a null subset of $\Hfngg$, the action of 
the automorphisms of $F_n$ preserve tori (products of the $\Uu(1)$ factors
and the invariant ellipses). Furthermore on almost every torus, the action
is ergodic, and hence every invariant function factors through the trace
functions of the generators. The rest of the proof is completely analogous.

\nocite{*}
\bibliographystyle{amsplain}
\bibliography{outfn}

\providecommand{\bysame}{\leavevmode\hbox to3em{\hrulefill}\thinspace}
\providecommand{\MR}{\relax\ifhmode\unskip\space\fi MR }
\providecommand{\MRhref}[2]{%
  \href{http://www.ams.org/mathscinet-getitem?mr=#1}{#2}
}
\providecommand{\href}[2]{#2}
\begin{thebibliography}{1}

\bibitem{BG}
R.~Benedetto and W.~Goldman, \emph{The topology of the relative character
  variety of the quadruply-punctured sphere}, Experimental Mathematics (1999),
  no.~8:1, 85--104.

\bibitem{Erg}
W.~Goldman, \emph{Ergodic theory on moduli spaces}, Ann. Math. (1997), no.~146,
  475--507.

\bibitem{}
\bysame, \emph{Action of the modular group on real
  $\operatorname{SL(2)}$-characters of a one-holed torus}, Geometry and
  Topology. (2003), no.~7, 443--486.

\bibitem{PickrellXia}
D.~Pickrell and E.~Xia, \emph{Ergodicity of mapping class group actions on
  representation varieties, I. Closed surfaces}, Comment.\
  Math.\ Helv.\ (2001), no.~77, 339--362.


\end{thebibliography}

%
%
%
%
%
%
%

\end{document}